\documentclass[11pt]{amsart}
\usepackage{amscd,amssymb,latexsym, verbatim}

\newtheorem{theorem}{Theorem}[section]

\theoremstyle{definition}
\newtheorem{definition}[theorem]{Definition}

\newtheorem{remark}[theorem]{Remark}
\newtheorem*{ack}{Acknowledgment}

\numberwithin{equation}{section}

\newcommand{\tsum}{\textstyle\sum\limits}

\newcommand{\A}{\mathcal{A}}

\newcommand{\RR}{\mathcal{R}}

\newcommand{\C}{\mathbb{C}}
\newcommand{\F}{\mathbb{F}}
\newcommand{\Z}{\mathbb{Z}}
\newcommand{\Q}{\mathbb{Q}}

\newcommand{\K}{\mathbb{K}}

\newcommand{\e}{{\epsilon}}
\newcommand{\g}{{\gamma}}

\renewcommand{\a}{{\alpha }}
\renewcommand{\b}{{\beta}}
\renewcommand{\d}{\delta}
\renewcommand{\l}{\lambda}

\newcommand{\md}[1]{\ (\bmod\:{#1})}

\begin{document}

\title[Massey Products of Complex Hypersurface Complements]%
{Massey Products of Complex Hypersurface Complements}

\author{Daniel~Matei}
\address{Graduate School of Mathematical Sciences, University of Tokyo, Tokyo, Japan
and Institute of Mathematics of the Romanian Academy, Bucharest, Romania.}



\begin{abstract}
It was shown by Kohno that all higher Massey products in the rational cohomology 
of a complex hypersurface complement vanish.
We show that in general there exist non-vanishing triple Massey products in the cohomology 
with finite field coefficients.
\end{abstract}

\maketitle
\section{Introduction}
\label{sec:intro}

The study of the topology of hypersurface complements is a classical subject in algebraic geometry.
Most of what is known about these spaces is related to invariants of their rational homotopy type. 
In this paper, we attempt to show that their $\F_p$-homotopy  type captures in general more information 
than the $\Q$-homotopy type, where $\F_p$ is the prime field of $p$ elements.
Let $X$ be the complement to a hypersurface $S$ in $\C P^{d}$. Then we have the following results due to
Kohno~\cite{Ko1, Ko2}: Massey products in $H^{*}(X,\Q)$ of length $\ge 3$ vanish. 
Moreover, the Malcev Lie algebra of $\pi_1(X)$ and the completed holonomy Lie algebra of $H^{\le 2}(X,\Q)$ are isomorphic. 
Thus, the $\Q$-completion of $\pi_1(X)$ is completely determined by the $\Q$-cohomology algebra of $X$. 
In the case when $S$ is a hyperplane arrangement $X$ is $\Q$-formal by Morgan~\cite{Mo},  
that is the entire $\Q$-homotopy type of $X$ is determined by the algebra $H^{*}(X,\Q)$. 
In this context, it seems natural to pose the following questions: Are there non-vanishing Massey products in $H^{*}(X,\F_p)$ 
for all primes $p$? Is $X$ a $\F_p$-formal space, particularly when $X$ is a hyperplane arrangement complement?

Massey products are known to be obstructions to formality, see~\cite{DGMS, El}. 
So, if the answer to the first question was yes, then the space $X$ would not be $\F_p$-formal. 
For compact K\"{a}hler manifolds the above questions were answered by Ekedhal in~\cite{Ek} by constructing 
such manifolds $M$ with non-vanishing triple products in $H^*(M, \F_p)$. 
Thus, a compact K\"{a}hler manifold although is $\Q$-formal by~\cite{DGMS}, in general it may not be $\F_p$-formal. 
The case of non-compact complex algebraic varieties is already different over $\Q$ from the compact case. 
As pointed out by Morgan in~\cite{Mo} such varieties may not be $\Q$-formal. 

The main result of this paper settles in affirmative the existence of
non-vanishing Massey products in the $\F_p$-cohomology of a hyperplane arrangement complement
for all odd primes $p$, thus showing that arrangement complements are not $\F_p$-formal
in general.
\begin{theorem}
\label{th:main}
For every odd prime $p$, the complement $X$ to the complex reflection arrangement $\A$ in $\C^3$ associated 
with the unitary reflection group $G(p,1,3)$ has, modulo indeterminacy, non-vanishing Massey products in $H^2(X,\F_p)$.
\end{theorem}
The cohomology operations that came to be known as Massey products were introduced by
W.~S.~Massey in~\cite{Mas}. Since then, they became important tools in algebraic topology, being especially
used as means of distinguishing spaces with the same cohomology but different homotopy type.
In general they are rather complicated objects, since they are in fact sets of cohomology classes.
But, in certain cases, they turn out to be cosets as shown by May in~\cite{May}, the simplest
instance being that of the Massey products of three cohomology classes.
In this paper we will only consider triple Massey products of cohomology classes of degree $1$ in 
the cohomology algebra in degrees at most $2$. 
In fact all the computations will take place in the group cohomology of the fundamental group 
$\pi_1(X)$ of our hypersurface complement. 
By the Lefschetz-Zariski classical theorem a generic $2$-dimensional section 
of $X$ captures all that topological information. 

The hypersurfaces $S$ that we will be our main focus are the hyperplane arrangements. 
Firstly because their complements are $\Q$-formal as discussed above.
Secondly because the integral cohomology of their complements is known to be 
torsion-free, see~\cite{OT}. In general, for a hypersurface $S$ consisting of non-linear
irreducible components, $H^{*}(X,\Z)$ will have torsion, and thus, at least conceptually, 
the chances of getting non-vanishing Massey products in $H^{*}(X,\F_p)$
are already greater. However, it is possible for a non-linear hypersurface to have 
torsion-free $H^{*}(X,\Z)$ as long as sufficiently many components of it are hyperplanes. 
We will briefly consider an example of such a non-linear hypersurface that nevertheless 
has triple non-vanishing Massey products. 

The arrangement $\A$ of Theorem~\ref{th:main} is the full monomial arrangement 
$\A(r, 1, 3)$ in $\C^3$, with $r=p\ge 3$.
These arrangements are members of a series of complex reflection arrangements, 
$\A(r, 1, d)$, associated to the full monomial reflection group  $G(r, 1, d)$, see~\cite{BMR, CS, OT}. 
For $r\ge 1$ and $d\ge 2$ let $\A(r, 1, d)$ be defined by the polynomial
$Q=z_1\cdots z_d \cdot\prod_{1\le i<j\le d} (z_i^r-z_j^r)$.
Note that the arrangements $\A(1, 1, d)$ and $\A(2, 1, d)$ are the Coxeter arrangements 
of type $A$ and respectively $B$.

For $\A$ a complex hyperplane arrangement in $\C^d$, with complement $X$ and group $G=\pi_1(X)$,
it is known that the rings $H^{*\le 2}(X,\K)$ and $H^{*\le 2}(G,\K)$ are isomorphic for 
$\K$ a field or $\Z$, see for example~\cite{MScn}.  Moreover, the complement to 
$\A(r, 1, d)$ is a $K(\pi, 1)$ with $\pi$ the pure braid group $P(r, 1, d)$ associated to $G(r, 1, d)$, 
see Orlik and Solomon~\cite{OS}. 
Taking advantage of this, we use the cochains of $G$ rather than those of $X$ to compute Massey products. 
We will use a presentation of $P(r, 1, d)$ obtained by Cohen~\cite{C}. 

A key r\^{o}le in the computations is played by the so-called resonance varieties of the arrangement,
see~\cite{Fa, MScn}. The resonance variety $\RR(\A,\K)$  over a field $\K$ of an arrangement $\A$ is 
the subvariety of $H^1(X,\K)$ encoding the vanishing cup products:
\[\RR(\A,\K)=\left\{ \l \in H^1(X,\K)\: \left|\:
\exists \,\mu \not\in \K\l \text{ such that } \l \cup \mu= 0 \right\}\right..
\]
The knowledge of the classes in $H^1(X,\K)$ that cup zero is especially needed for calculating 
a triple Massey product $\langle\a, \b, \g\rangle$ as that is well-defined only when $\a\cup\b=\b\cup\g=0$. 
In~\cite{Fa}, Falk gives a combinatorial recipe to detect posible essential components of  $\RR(\A,\K)$.
For $\A=\A(r, 1, 3)$, the classes used to define the non-vanishing Massey products belong to 
such components arising when $\K=\F_p$, for the special primes $p$ dividing $r$.
It can be shown that $\A(r, 1, 3)$ presents non-vanishing $\F_p$-Massey products for all 
primes $p$ and all multiples $r$ of $p$ (multiples of $4$ if $p=2$). Here only the case $r=p$ is treated.

The paper is organized as follows. In Section~\ref{sec:massprod} we define the triple Massey products of 
a $2$-complex associated to a finitely presented group, and explain how they can be computed from the presentation. 
In Section~\ref{sec:reflarr} we introduce the monomial arrangements and give presentations by generators and 
relators of the fundamental groups of their complements. In Section~\ref{sec:nvanmp} we exhibit non-vanishing
triple Massey products in the $\F_p$-cohomology of the complements to
$3$-dimensional monomial arrangements, for $p$ an odd prime. We also present 
a non-linear arrangement of curves in $\C P^2$ whose complement has 
non-vanishing triple Massey products over $\Z_2$. In the last section we pose
some further questions that we intend to explore elsewhere.

\section{Massey products of $CW$-complexes} 
\label{sec:massprod}
The results on Massey products that we need may be found in the works of Porter~\cite{Po}, Turaev~\cite{Tu}, 
and Fenn and Sjerve~\cite{FS1, FS2}. In these papers the Massey products of $1$-cohomology classes are computed in terms of the 
so-called Magnus coefficients, via the free calculus of Fox.
Unless otherwise specified, all the homology and cohomology groups 
will have coefficients in $\F_p$, the integers modulo a prime $p$.

\begin{definition}
\label{def:masprod}
Let $X$ be a space of the homotopy type of a $CW$-complex. If $\a, \b, \g$ in $H^1(X)$ are such that 
$\a\cup \b=\b\cup \g=0$ then the triple Massey product $\langle \a, \b, \g\rangle$ is defined as follows:
Choose representative $1$-cocycles $\a', \b', \g'$ and 
cochains $x, y$ in $C^1(X)$ such that $dx=\a'\cup \b'$ and $dy=\b'\cup \g'$. 
Then $z=\a'\cup y+x\cup \g'$ is a $2$-cocycle. The cohomology classes $z\in H^2(X)$ constructed in this way
are only determined up to $\a\cup H^1(X) + H^1(X)\cup \g$, and they form a set denoted by  
$\langle \a, \b, \g\rangle$.  

As pointed out by May~\cite{May}, the indeterminacy  is a vector space, and so $\langle \a, \b, \g\rangle$ can be 
thought of as a coset modulo $\a\cup H^1(X) + H^1(X)\cup \g$. 
The triple Massey product $\langle \a, \b, \g\rangle$ is said to be vanishing if this coset is trivial.
\end{definition} 

In this paper $X$ will always be a $K(G, 1)$ for  $G$ a finitely presented group. 
We will identify from now on the cohomology of $X$ with that of $G$.

Let $G=\langle x_1, \dots, x_n \mid R_1, \dots, R_{m}\rangle$ be a presentation for $G=\pi_1(X)$.
Assume that $R_l$ is a commutator and that the presentation is minimal. 
By Hopf's formula the homology classes of the relators $R_l$ form a basis in $H_2(G,\Z)=\Z^{m}$.
Morover, the generators $x_i$ determine a basis of $H_1(G,\Z)=\Z^n$. 
Let $e_i$ be the dual basis in $H^1(G,\Z)=\Z^n$. 

Let $F$ be the free group on $x_1, \dots, x_n$. 
If $w$ is a word in $F$ then its Fox derivative $\partial_j(w)$ is computed by the following rules:
$\partial_j (1)=0, \,\partial_j (x_i)=\delta_{i, j}$, and $\partial _j (uv)=\partial_j (u)\epsilon (v)+u\partial_j (v)$,
where $\e:\Z F\to\Z$ is the augmentation of the group ring $\Z F$.

Let $I=(i_1,\dots, i_q)$ be a multi-index with $i_j$ taking values in $1,\dots, n$.
The Magnus $I$-coefficient of a word $w$ is defined by $\e_{I}^{(0)}(w)=\e\partial_{I}(w)$, where
$\partial_{I}(w)=\partial_{i_1}\dots\partial_{i_q}(w)$. The  $\F_p$-valued Magnus coefficients $\e^{(p)}_I(w)$ 
of $w$ are defined simply by taking integers modulo the prime $p$. Most of the time
we will drop the reference to it and simply write $\e_{I}(w)$ for $\e^{(p)}_I(w)$.
We will usually refer to $\epsilon_{i, j}(w)$ and
$\epsilon_{i, j, k}(w)$ as a double, and respectively triple Magnus coefficient.

The following result~\cite{FS1, FS2, Po, Tu} will be used to compute triple Massey products.
Let $\a, \b, \g$ be cohomology classes in $H^1(G)$ such that $\a\cup \b=\b\cup \g=0$. 
Denote by $(\cdot,\cdot)$ the Kronecker pairing between cohomology and homology. Then we have:
\begin{theorem}
The Massey product $\langle \a, \b, \g\rangle$ of $\a=\sum \a_i e_i$, $\b=\sum \b_j e_j$, 
$\g=\sum \g_k e_k$ contains $\xi$, where:
\[
(\xi, R_l)=\sum_{1\le i, j, k \le n} \a_i \b_j \g_k \cdot\epsilon_{i, j, k}(R_l).
\]
\end{theorem}
From now on by $\langle \a, \b, \g\rangle$ we will understand the coset of the class 
$\xi$ modulo the indeterminacy $\a\cup H^1(G) + H^1(G)\cup \g$. 
The Massey product $\langle \a, \b, \g\rangle=\xi$ is functorial with respect to maps of spaces, 
as shown by Fenn and Sjerve in~\cite{FS1, FS2}.

The following formulae can readily deduced from the definitions and they will 
be used to compute  the Magnus coefficients of a commutator word.
\begin{equation}
\label{eq:e2com}
\e_{k, l}([u, v])=\e_{k}(u)\e_{l}(v)-\e_k(v)\e_{l}(u), \text{ where } [u, v]=uvu^{-1}v^{-1}.
\end{equation}
\begin{align}
\label{eq:e3com}
\e_{k, l, m}([u, v])=
&\e_{k}(u)\e_{l, m}(v)-\e_{m}(u)\e_{k, l}(v)+\e_{k, l}(u)\e_{m}(v)-\e_k(v)\e_{l, m}(u)+\\
\notag
&\left(\e_{k}(v)\e_l(u)-\e_{k}(u)\e_{l}(v)\right)\cdot\left(\e_{m}(u)+\e_{m}(v)\right).
\end{align}
We will also need formulae for products of conjugated generators:
\begin{equation}
\label{eq:e1cnj}
\e_{k}(x_{i_1}^{w_1}\cdots x_{i_j}^{w_j})=\sum_{a=1}^{j}\e_{k}(x_{i_a})=\sum_{a=1}^{j}\delta_{k, i_a}.
\end{equation}
\begin{equation}
\label{eq:e2cnj}
\e_{k, l}(x_{i_1}^{w_1}\cdots x_{i_j}^{w_j})=
\sum_{a=1}^{j}\left(\e_{k}(w_{a})\delta_{l, i_a}-\e_{l}(w_{a})\delta_{k, i_a}\right)+
\sum_{1\le a < b\le j}\delta_{k, i_a}\delta_{l, i_b},
\end{equation}
where $x^a=axa^{-1}$ and $\delta_{i,j}$ is Kronecker's delta.

\section{Monomial arrangements and their groups}
\label{sec:reflarr}

We introduce in this section our main examples of hypersurfaces whose complements 
have non-vanishing Massey products in the $\F_p$-cohomology. 
They are the complex reflection 
arrangements $\A(r, 1, d)$ associated with the monomial reflection group 
$G(r, 1, d)$. Their complements are $K(\pi, 1)$ spaces and for all practical purposes we will 
identify their cohomology with that of their fundamental groups. We will describe here the group 
presentations that will be used in the Massey products computation.

\subsection{Arrangements groups}
\label{subsec:arrgrp}
We start with a brief
overview of the fundamental group of hyperplane complements.
Let $\A$ be a hyperplane arrangement in the affine space $\C^d$ and $X$ its complement.
Let us recall now the most salient features of the fundamental group $G=\pi_1(X)$ as a finitely presentable 
group. For all the details see~\cite{OT}. 
First $G$ is generated by the meridians $\g_H$ around each hyperplane $H\in \A$.
Each codim $2$ 
intersection $H_{i_1}\cap\dots\cap H_{i_n}$ of  hyperplanes in $\A$ 
determines $n-1$ relations: 
$g_1g_2 \cdots g_n=g_2 \cdots g_n\cdot g_1=\cdots=g_n\cdot g_1\cdots g_{n-1}$, 
where $ g_j$ is some conjugate of the generator $x_j=g_{H_j}$. 
We denote  by $[g_1,\dots, g_n]$
the family of commutator relators $[g_1 \dots g_i, g_{i+1} \dots g_n]$, 
with $1\le i <n$. 

Thus we are lead to compute the Magnus coefficients of relators in families of the form:
$\big[x_{i_1}^{w_1},\dots, x_{i_n}^{w_n}\big]$. Note that the indices $i_j$ are all distinct.
Denote by $R_{I, w}^{j}$ the commutator 
$\big[x_{i_1}^{w_1}\dots x_{i_j}^{w_j}, x_{i_{j+1}}^{w_{j+1}}\dots, x_{i_n}^{w_n}\big]$.

The Magnus coefficients of order $2$ of $R_{I, w}^{j}$ are given by:
\begin{equation}
\label{eq:epstwo}
\e_{k, l}\left(R_{I, w}^{j}\right)=\sum_{a=1}^{j}\sum_{b=j+1}^{n} 
(\d_{k, i_{a}}\d_{l, i_{b}}-\d_{k, i_{a}}\d_{k, i_{b}}).
\end{equation}
It is easily seen that:
\begin{equation}
\notag
\e_{k, l}\left(R_{I, w}^{j}\right)=
\begin{cases}
1  &\text{if $k=i_a$ and $l=i_b$, for some $1\le a\le j$ and $j+1\le b\le n$}\\
-1  &\text{if $k=i_b$ and $l=i_a$, for some $1\le a\le j$ and $j+1\le b\le n$}\\
0   &\text{otherwise}. 
\end{cases}
\end{equation}

The Magnus coefficients of order $3$ of $R_{I, w}^{j}$ are given by:
\begin{equation}
\label{eq:epsthree}
\e_{k, l, m}\left(R_{I, w}^{j}\right)=
\begin{cases}
\e_m(w_b)  &\text{if $k=i_a, l=i_b$, $m\not\in I$}\\
-\e_m(w_a)  &\text{if $k=i_b, l=i_a$, $m\not\in I$}\\
\e_l(w_b)+\e_k(w_{a'})-\e_l(w_a)+\delta_{a\le a'}  &\text{if $k=i_a$, $l=i_{a'}$, $m=i_b$}\\
\e_l(w_a)+\e_k(w_{b'})-\e_l(w_b)+\delta_{b\le b'}  &\text{if $k=i_b$, $l=i_{b'}$,  $m=i_a$}\\
\e_k(w_b)-\e_m(w_b)-1 &\text{if $k=i_a$, $l=i_b$, $m=i_{a'}$}\\
\e_l(w_{b'})-\e_l(w_a)-\e_m(w_b)+\delta_{b\le b'}-1  &\text{if $k=i_a$, $l=i_b$, $m=i_{b'}$}\\
\e_k(w_a)-\e_m(w_a)+1 &\text{if $k=i_b$, $l=i_a$, $m=i_{b'}$}\\
\e_l(w_{a'})-\e_l(w_b)-\e_m(w_a)+\delta_{a\le a'}+1  &\text{if $k=i_b$, $l=i_a$, $m=i_{a'}$}\\
0   &\text{otherwise}.
\end{cases}
\end{equation}
where always $1\le a, a'\le j$ and $j+1\le b, b'\le n$.

\subsection{Monomial arrangements}
\label{subsec:monarr}
We introduce now our main class of examples. 
For $r\ge 1$ and $d\ge 2$ let  $\A(r, 1, d)$ be the arrangement defined by:
\[
Q=z_1\cdots z_d\cdot\prod_{1\le i<j\le d} (z_i^r-z_j^r). 
\]
The complement of $\A(r, 1, d)$ is a $K(\pi, 1)$ with $\pi$ the pure braid group 
$P(r, 1, d)$ associated to the full monomial complex reflection group  $G(r, 1, d)$, 
see~\cite{BMR, OS}. The group $P(r, 1, d)$ admits an iterated semidirect product 
structure:
$P(r, 1, d)=F_{n_l}\rtimes\cdots\rtimes F_{n_1}$, where $n_i=(i-1)r+1$ for $1\le i\le r$, 
as shown in~\cite{BMR}.

A presentation for $P(r, 1, d)$ was obtained by Cohen in~\cite{C}.
Following that paper, let us first describe the codim $2$ intersections among the hyperplanes of 
$\A(r, 1, d)$: 
\begin{align}
&H_i\cap H_{i,j}^{(1)}\cap\cdots\cap H_{i,j}^{(r-1)}\cap H_j\cap H_{i,j}^{(r)}\\
&H_k\cap H_{i,j}^{(q)}, \text{if $k\neq i$ or $k> j$}\\
&H_{i,j}^{(q)}\cap H_{k,l}^{(s)}, \text{if $i, j, k, l$ distinct}, \\
&H_{i,j}^{(q)}\cap H_{j,k}^{(s)}\cap H_{i,k}^{(t)}, \text{ if $t=q+s \md{r}$}, 
\end{align}
where $H_i=\{z_i=0\}, 1\le i\le 3$, and $H_{i,j}^{(q)}=\{z_i=\zeta^q z_j\}$, 
where $\zeta=\exp(2\pi i/r)$, $1\le i< j\le 3$ and $1\le q\le r$.

We focus now on the case $d=3$. In~\cite{C} a presentation of $P(r, 1, 3)$ is given having
$3r+3$ generators, say $x_1,\dots, x_{3r+3}$, and $2r^2+6r+3$ relators.
We group these relators in nine families corresponding to the types of the codimension $2$ intersections.
\begin{align}
\label{eq:relsfirst}
&A=\big[x_{3r+1}, x_1,\dots, x_{r-1}, x_{3r+2}, x_r\big], \\
&B=\big[x_{3r+1}, x_{2r+1},\dots, x_{3r-1}, x_{3r+3}, x_{3r}\big], \\
&C=\Big[x_{3r+2}, x_{r+1}^{x_r x_{3r+1}x_1 x_2 \cdots x_{r-1}x_{3r+1}^{-1}},\dots, 
x_{2r-1}^{x_r x_{3r+1}x_1 x_{3r+1}^{-1}}, x_{3r+3}, x_{2r}\Big], \\
&D_{1, s}=\big[x_{3r+1}, x_{r+i}\big], 1\le s\le r,\\
&D_{2, s}=\big[x_{3r+3}, x_{i}\big], 1\le s\le r,\\
&D_{3, s}=\Big[x_{3r+2}^{x_i x_{i+1} \cdots x_{r-1}}, x_{2r+i}^{x_{2r}}\Big], 1\le s\le r,\\
&T_{s}=\big[x_{s}, x_{2r+s}, x_{2r}\big], 1\le s\le r, \\
&U_{t,s}=\Big[x_{s}, x_{2r-t}, x_{2r+s-t}^{x_{2r-t+1} \cdots x_{2r-1}}\Big], 1\le t < s\le r,  \\
\label{eq:relslast}
&V_{s,t}=\Big[x_{s}^{x_{3r+1}}, x_{2r-t}, x_{3r+s-t}^{x_{2r-t+1} \cdots x_{2r-1}}\Big], 1\le s \le t < r.
\end{align}
Now, recall that the notation $R=[x_{i_1}^{w_1},\dots, x_{i_n}^{w_n}]$ stands for the following set of commutators:
$\{R^j=[x_{i_1}^{w_1}\dots x_{i_j}^{w_j}, x_{i_{j+1}}^{w_{j+1}}\dots, x_{i_n}^{w_n}] 
\mid 1\le j < n\}$. Thus,  in agreement with the notations of~\eqref{eq:relsfirst}-\eqref{eq:relslast},
the relators in $P(r, 1, 3)$ will be denoted by: 
$A^j, B^j, C^j$, where $j=1, \dots, r+1$, and
$D_{1, s}, D_{2, s}, D_{3, s}$, where $j=1$ and is omitted,
and finally $T_s^j, U_{t,s}^j, V_{s,t}^j$, where $j=1, 2$.

\section{Non-vanishing triple Massey products}
\label{sec:nvanmp}

In this section we will present non-vanishing triple Massey products in the $\F_p$ cohomology of certain 
hypersurface complements. All such products will be of the form $\langle \a, \a, \b \rangle$ with $\a$ 
and $\b$ linearly independent. The main example will consist of the monomial arrangements introduced in 
the previous section. We will also give an example of a non-linear arrangement 
of curves with the desired non-vanishing property.

\subsection{Resonance varieties} 
\label{subsec:resvar}
We first determine the vanishing cup products in $H^2(X, \F_p)$, for $X$ the 
complement of
a monomial arrangement $\A$, using an invariant of a cohomology ring introduced by Falk in~\cite{Fa}.
The resonance variety $\RR(\A,\F_p)$ of an arrangement $\A$ is the subvariety of $H^1(X,\F_p)$ defined by:
\[
\RR(\A,\F_p)=\left\{ \l \in H^1(X,\K)\: \left|\:
\exists \,\mu \not\in \K\l \,\text{ such that } \,\l \cup \mu= 0 \right\}\right..
\]
In~\cite{Fa} it is shown how one can construct components of $\RR(\A, \F_p)$ from 
the so-called neighborly partitions of the arrangement $\A$. The neighborly partitions of the monomial 
arrangements $\A=\A(r, 1, 3)$ have been determined in~\cite{CS}. The most interesting for us is the partition 
$\Pi=(H_3, H_{12}^{(i)}\mid H_2, H_{13}^{(j)}\mid H_1, H_{23}^{(k)})$ giving
rise to a component 
$C_{\Pi}$ of $\RR(\A, \F_p)$ having the following equations:
\begin{align}
\label{eq:cpi}
\notag
&\l_1+\dots+\l_r=\l_{r+1}+\dots+\l_{2r-1}=\l_{2r+1}+\dots+\l_{3r-1}=0 \\  
\notag
&\l_i+\l_{2r}+\l_{2r+i}=0, \quad 1\le i\le r, \\ 
&\l_i+\l_{2r-j}+\l_{2r+i-j}=0, \quad 1\le j < i\le r, \\ \notag
&\l_i+\l_{2r-j}+\l_{3r+i-j}=0, \quad 1\le i\le j < r, \\ \notag
&\l_{3r+1}=\l_{3r+2}=\l_{3r+3}=0
\end{align}
It is easily seen that $\dim C_{\Pi}=3$ if $p$ divides $r$ (or $4$ divides $r$, if $p=2$), and $\dim C_{\Pi}=2$, otherwise.

\subsection{Massey products of monomial arrangements} 
\label{subsec:mpmonarr}
We prove here the main result, showing that, in general, Massey products in the positive characteristic cohomology of a 
hypersurface complement may not vanish modulo indeterminacy, although over the rationals they always do so.

\begin{theorem}
For every odd prime $p$ the complement $X$ of the arrangement $\A(p, 1, 3)$ 
in $\C^3$ of degree $3p+3$ has non-vanishing triple Massey products in $H^2(X, \F_p)$.
\end{theorem}

\begin{proof}
We will show that a certain triple product $\langle \a, \a, \b \rangle$ does not vanish modulo its indeterminacy.
The cohomology classes $\a$ and $\b$ are given in coordinates by
\[
a:\:\a_i=1, \a_{r+i}=-1, \a_{2r+i}=\a_{3r+1}=\a_{3r+2}=\a_{3r+3}=0,
\]
and respectively by 
\[
\b:\:\b_i=0, \b_{r+i}=1, \b_{2r+i}=-1, \b_{3r+1}=\b_{3r+2}=\b_{3r+3}=0,
\]
where $1\le i\le r$. Clearly the points $\a$ and $\b$ satisfy the equations~\eqref{eq:cpi}, so they belong 
to $C_{\Pi}$, and moreover $\a\cup\b=0$.
Using~\eqref{eq:epsthree} we can express $\langle \a, \a, \b \rangle$ in the basis
of $H^2(X, \F_p)$ given by the duals of the relators~\eqref{eq:relsfirst}-\eqref{eq:relslast}, 
abusing the notation for the sake of simplicity.

\begin{align}
\label{eq:aab}
\langle \a, \a, \b \rangle=
&\sum_{j=1}^{p} (j-1) C^{j}+(p-1) C^{p+1}-\sum_{s=1}^{p}T_s^2+
\sum_{1\le t<s\le p} t\, U_{t,s}^1+\sum_{1\le t<s\le p} U_{t,s}^2+\\ \notag
&\sum_{1\le s\le t< p} t\, V_{s,t}^1+\sum_{1\le s\le t< p} V_{s,t}^2. 
\end{align}

Next, using~\eqref{eq:epstwo}, we obtain the indeterminacy 
$\a\cup H^1(X) + H^1(X)\cup \b$. 
If $a=\tsum a_i e_i$ and $b=\tsum b_i e_i$ are arbitrary classes
in $H^1(X)$ then we find the following expression for $\a\cup a + b\cup\b$:

\small
\begin{align}
\label{eq:ind}\notag
&\qquad\Bigg(\sum_{s=1}^{p} a_s+a_{3p+1}+a_{3p+2}\Bigg)
\Bigg(\sum_{j=1}^{p} (j-1) A^{j}+(p-1) A^{p+1}\Bigg)+
\Bigg(\sum_{s=1}^{p} b_{2p+s}+b_{3p+1}+b_{3p+3}\Bigg)\\ \notag
&\qquad\Bigg(\sum_{j=1}^{p} (j-1) B^{j}+(p-1) B^{p+1}\Bigg)-
\Bigg(\sum_{s=1}^{p} (a_{p+s}+b_{p+s})+
a_{3p+2}+a_{3p+3}+b_{3p+2}+b_{3p+3}\Bigg) \\ \notag
&\qquad\Bigg(\sum_{j=1}^{p} (j-1) C^{j}+(p-1) C^{p+1}\Bigg)+
\left(a_{3p+1}+b_{3p+1}\right)\sum_{s=1}^{p} D_{1,s}-a_{3p+3} 
\sum_{s=1}^{p} D_{2,s}-\\ 
&\qquad b_{3p+2}\sum_{s=1}^{p} D_{3,s}+
\sum_{s=1}^{p}\left(a_s+a_{2p+s}+a_{2p}\right)\left(T_s^1+T_s^2\right)+
\sum_{s=1}^{p}\left(b_s+b_{2p+s}+b_{2p}\right) T_s^2 +\\ \notag
&\qquad\sum_{1\le t<s\le p}\left(a_s+a_{2p-t}+a_{2p+s-t}\right) U_{t,s}^1-
\sum_{1\le t<s\le p}\left(b_s+b_{2p-t}+b_{2p+s-t}\right) U_{t,s}^2+\\ \notag
&\qquad\sum_{1\le s\le t< p}\left(a_s+a_{2p-t}+a_{3p+s-t}\right) V_{s,t}^1-
\sum_{1\le s\le t< p}\left(b_s+b_{2p-t}+b_{3p+s-t}\right) V_{s,t}^2. 
\end{align}
\normalsize
We want to show that the triple Massey product $\langle \a, \a, \b \rangle$ 
does not vanish modulo indeterminacy. Suppose that it does vanish, and so 
there exist $a$ and $b$ in $H^1(X)$ such that $\langle \a, \a, \b \rangle$ is of 
the form $\a\cup a + b\cup\b$. This leads to the following set of equations over $\F_p$:
\begin{align}
\label{eq:abeq1}
&\sum_{s=1}^{p} a_s+a_{3p+1}+a_{3p+2}=\sum_{s=1}^{p} b_{2p+s}+b_{3p+1}+b_{3p+3}=0, \\
\label{eq:abeq2}
&\sum_{s=1}^{p} (a_{p+s}+b_{p+s})+a_{3p+2}+a_{3p+3}+b_{3p+2}+b_{3p+3}=-1, \\
\label{eq:abeq3}
&a_{3p+1}+b_{3p+1}=a_{3p+3}=b_{3p+2}=0, \\
\label{eq:abeq4}
&a_s+a_{2p+s}+a_{2p}=0,\,  b_s+b_{2p+s}+b_{2p}=-1,\\
\label{eq:abeq5}
&a_s+a_{2p-t}+a_{2p+s-t}=t,\, b_s+b_{2p-t}+b_{2p+s-t}=-1, \\
\label{eq:abeq6}
&a_s+a_{2p-t}+a_{3p+s-t}=t,\, b_s+b_{2p-t}+b_{3p+s-t}=-1, 
\end{align}
where the ranges of the indices in~\eqref{eq:abeq4},~\eqref{eq:abeq5}, and~\eqref{eq:abeq6} are those 
in~\eqref{eq:ind}.

Now from~\eqref{eq:abeq4},~\eqref{eq:abeq5}, and~\eqref{eq:abeq6} we can readily see that we 
must have: 
\[
\sum_{s=1}^{p} a_s=\sum_{s=1}^{p} a_{p+s}=\sum_{s=1}^{p} a_{2p+s}=0
\text{ and }
\sum_{s=1}^{p} b_s=\sum_{s=1}^{p} b_{p+s}=\sum_{s=1}^{p} b_{2p+s}=0.
\]
From these equations combined with \eqref{eq:abeq1}, \eqref{eq:abeq2}
and \eqref{eq:abeq3} we obtain:
$a_{3p+1}+a_{3p+2}=b_{3p+1}+b_{3p+3}=a_{3p+1}+b_{3p+1}=0$ and $a_{3p+2}+b_{3p+3}=-1$.
But this system of equations clearly has no solution.

\end{proof}

\begin{remark}
We will show elsewhere that in fact any triple Massey product in $H^2(P,\F_p)$ of the form 
$\langle \a, \a, \b \rangle$ with $\a$ and $\b$ (not proportional) in 
$C_{\Pi}\subset \RR(P, \F_p)\subset H^1(P,\F_p)$ does not vanish modulo 
the indeterminacy $\a\cup H^1(P,\F_p)+H^1(P,\F_p)\cup\b$, if $p\mid r$ 
(or $4\mid r$, if $p=2$), where $P=P(r, 1, 3)$. Thus it will follow that for 
every prime $p$ and multiple $N\ge 3$ of $p$ 
(of $4$ if $p=2$) there exists a line arrangement $\A$ in $\C^2$ of degree   
$3N+3$ whose complement $X$ has non-vanishing triple Massey products in $H^2(X, \F_p)$.
\end{remark}

\subsection{Curves with non-linear components}
\label{subsec:curves}

Let $\mathcal{C}=Q_2\cup T_1\cup T_2\cup T_3$ be the curve in $\C P^2$ of degree $5$, consisting
of a smooth irreducible curve $Q_2$ of degree $2$ and three lines $T_1, T_2, T_3$ tangent to $Q_2$. 
As explained by Kaneko, Tokunaga and Yoshida in~\cite{KTY},
this curve is related with the discriminant 
of a certain crystallographic group, thus is of the same nature as the above reflection arrangements. 
In~\cite{KTY} a presentation for the fundamental group of the complement to $\mathcal{C}$ in $\C P^2$ is determined:
\[
\pi_1(\C P^2\setminus\mathcal{C})=\langle x_1, x_2, x_3 \mid [x_3x_ix_3,x_i], i=1,2, [x_3x_1x_3^{-1},x_2] \rangle.
\]
An easy computation with double Magnus coefficients shows that all $\Z_2$ cup products 
$e_i\cup e_j$ vanish except for $e_1\cup e_2$.
Moreover, by computing triple Magnus 
coefficients we can see that the Massey products $\langle \a, \a, \b \rangle$ over $\Z_2$ 
do not vanish, if $\a\not\in\Z_2\cdot(e_1+e_2+e_3)$.

\begin{remark}
It is possible to generalize this example to a curve 
$\mathcal{C}=Q_d\cup T_1\cup\cdots\cup T_n$ of degree $d+n$, where $Q_d$ is a smooth irreducible 
curve of degree $d\ge 2$ and $T_1\dots T_n$ are $n\ge d+1$ tangent lines to $Q$. Then the complement $X$ of 
$\mathcal{C}$ will have non-vanishing triple Massey products of the form $\langle \a, \a, \b \rangle$ in 
$H^2(X, \F_p)$, for every prime $p$ dividing $d$.
\end{remark}

\section{Further questions}
\label{sec:quest}
We end the paper by raising a few questions:

\begin{enumerate}
\item Note that the above arrangements exhibiting non-vanishing Massey products do not admit linear equations over the reals! 
Is it true that real complexified arrangements never have non-vanishing Massey products? Computational evidence suggest that in 
this case Massey products in $H^2(X,\F_p)$ indeed all vanish.
\item All non-orientable matroids realizable over some $\Q(\alpha)$ lead to complex arrangements with 
non-vanishing Massey products?
\item Is there an analogue of Kohno's result over $\F_p$? Is it true that non-vanishing of higher Massey 
products over $\F_p$ 
implies that the $\F_p$-completion of $\pi_1(X)$ is not isomporhic to the completed holonomy algebra of $H^{\le 2}(X,\F_p)$?
\item Do non-linear curves (with enough cohomology) always present non-vanishing Massey products?

\item Are there any good criteria for $\F_p$-formality of $X$?  In this context, what is the r\^{o}le played by  
the lower $p$-central series of $\pi_1(X)$?
 
\end{enumerate}

\medskip
\begin{ack}
Most of this work was done while the author was visiting the University of Tokyo during the academic year $2003$.
He is grateful for all the support the Graduate School of Mathematical Sciences
has given to him during that time, in particular for making possible his participation
at the 2003 MSJ Conference ``Singularity theory and its applications'' in Sapporo,
where this material was first presented.

The computations in this paper were carried out with the help of the software 
package {\sl  Mathematica~5.0}.
\end{ack}

\end{document}